\documentclass[a4paper]{article}
\def\today{30.1.14} %PLEASE PUT HERE THE DATE
\usepackage{amsmath,amsfonts,amsthm,amssymb,amscd}
\usepackage{a4wide}
\theoremstyle{plain} \newtheorem{theorem}{Theorem}[section]
\newtheorem{lemma}[theorem]{Lemma}

 \theoremstyle{definition}
\newtheorem{definition}[theorem]{Definition} \theoremstyle{remark}
\newtheorem{remark}[theorem]{Remark}

\newcommand{\R}{{\mathbb R}} 
\newcommand{\Z}{{\mathbb Z}}
\newcommand{\Na}{{\mathbb N}}

\newcommand{\C}{\mathbb{C}}

\newcommand{\B}{{\mathcal B}}

\newcommand{\G}{{\mathcal G}}

\newcommand{\Ph}{{\mathcal P}}

\newcommand{\resto}{{\mathcal R}}

\newcommand{\Sc}{{\mathcal S}}
\newcommand{\U}{{\mathcal U}}

\newcommand{\tond}[1]{{\left(#1\right)}}
\newcommand{\quadr}[1]{{\left[#1\right]}}

\def\catene{BamG93}
\def\tedeschi{GauHL12}

\def\im{{\rm i}}

\def\uno{{\kern+.3em {\rm 1} \kern -.22em {\rm l}}}

\def\di{{\rm d}}
\def\norma#1{\left\| #1\right\|}

\def\poisson#1#2{\left\{#1;#2\right\}}

\def\tq{{\tt q}}
\def\tp{{\tt p}}
\def\tQ{{\tt Q}}

\def\normati#1{\left\| #1\right\|_{\sim}}

\numberwithin{equation}{section}

\begin{document}

\title{Normal Form and Energy Conservation \\ of High Frequency Subsystems
\\ Without Nonresonance Conditions}
\author{Dario Bambusi\footnote{Dipartimento di Matematica,
    Universit\`a degli studi di Milano, Via Saldini 50, 20133 Milano, Italy} \and Antonio Giorgilli$^{*}$ \and Simone Paleari$^{*}$ \and Tiziano Penati$^{*}$}

\date{\today}
\maketitle

%\listoftodos

%%%%%%%%%%%%%%%%%%%%%%%%%%%%%%%%%%%%%%%%%%%%%%%%%%%%%%%%%%%%%%%%%%%%%%%%
%
%               ABSTRACT
%
%%%%%%%%%%%%%%%%%%%%%%%%%%%%%%%%%%%%%%%%%%%%%%%%%%%%%%%%%%%%%%%%%%%%%%%%

\begin{abstract}
We consider a system in which some high frequency harmonic oscillators
are coupled with a slow system.  We prove that up to very long times
the energy of the high frequency system changes only by a small
amount. The result we obtain is completely independent of the
resonance relations among the frequencies of the fast system. More in
detail, denote by $\epsilon^{-1}$ the smallest high frequency. In the
first part of the paper we apply the main result of \cite{\catene} to
prove almost conservation of the energy of the high frequency system
over times exponentially long with $ {\epsilon^{-1/n}} $ ($n$ being
the number of fast oscillators). In the second part of the paper we
give e new self-contained proof of a similar result which however is
valid only over times of order $\epsilon^{-N}$ with an arbitrary
$N$. Such a second result is very similar to the main result of the
paper \cite{\tedeschi}, which actually was the paper which stimulated
our work.
\end{abstract}

%%%%%%%%%%%%%%%%%%%%%%%%%%%%%%%%%%%%%%%%%%%%%%%%%%%%%%%%%%%%%%%%%%%%%%%%
%
%               INTRODUCTION
%
%%%%%%%%%%%%%%%%%%%%%%%%%%%%%%%%%%%%%%%%%%%%%%%%%%%%%%%%%%%%%%%%%%%%%%%%

\section{Introduction}\label{intro}

In the phase space $\R^{2n}\oplus\R^{2d}\ni((p,q),(P,Q))$ we consider a
Hamiltonian system of the form
\begin{equation}
\label{sys.1}
H(p,q,P,Q)=h_\omega(p,q)+H_0(P,Q,q)\ ,
\end{equation}
where 
\begin{equation}
\label{homega.1}
h_{\omega}(p,q):=\sum_{j=1}^{n}\frac{p_j^2+\omega_j^2q_j^2}{2}
\end{equation}
is a system of ``fast'' harmonic oscillators and $H_0$ is an
analytic function describing a ``slow'' system (with canonical
variables $P,Q$) and its interaction with the fast system. We are
interested in the case where the frequencies $\omega_j$ are large, so
we define
\begin{equation}
\label{eps.1}
\epsilon:=\frac{1}{\min_{j}\{\omega_j\}}\ ,
\end{equation}
and study the system in the limit $\epsilon\to0$. 

In the first part of paper we apply the main result of \cite{\catene}
to prove that $h_\omega$ changes by a quantity which is at most of
order $\epsilon^{1/n}$ up to times exponentially long with
$\epsilon^{-1/n}$; in the second part we give a new self contained
proof of a stability result very close to a result by Gauckler, Heirer
and Lubich \cite{\tedeschi}, which ensures almost invariance of
$h_\omega$ over times of order $\epsilon^{-N}$ with an arbitrary
$N$. The main point is that all the results are completely
independent of the resonance relations among the frequencies $\omega_j$,
and thus hold uniformly for all the frequency vectors outside a cube
of side $\epsilon^{-1}$.

We recall that systems of the kind of \eqref{sys.1} arise in many
contexts; here we just mention the problem of the realization of
Holonomic constraints, in which the constraints are modeled by very
hard springs and one is interested in controlling if the dynamics of
the slow system converges, as $\epsilon\to0$ to the dynamics of
$H_0(P,Q,0)$. This is a very subtle question and indeed it is well
known that, in general, the convergence of the orbits can occur only
for times of order $\epsilon^{-1}$. For longer time scales one can
only pursue weaker results, and actually in \cite{BenGG87,BenGG89}
(see also \cite{BamG93}), it has been show that if the frequencies
$\omega_j$ are either completely resonant or fulfill some Diophantine
type inequalities, then $h_\omega$ is an approximate integral of
motion for times exponentially long with $\epsilon^{-a}$ with $a$
depending on the resonance properties of the frequency vector
$\omega$. All the constant involved in the main theorems of
\cite{BenGG87,BenGG89} depend on the properties of good/bad
approximability of the frequencies by rational vectors.

In their paper \cite{\tedeschi} Gauckler, Heirer and Lubich used
multiscale expansion to show that by restricting attention to time
scales of order $\epsilon^{-N}$ with arbitrary $N$, one can find
a result independent of the resonance properties of the frequencies,
and thus uniform for all frequencies outside an $n$ dimensional
hypercube of side $\epsilon^{-1}$. The paper \cite{\tedeschi} was
actually what stimulated the present work.
\vskip 10pt

In the present paper we present two results.

In the first part (Section~\ref{1}) we look for stability over
exponentially long times, in the spirit of Nekhoroshev theory. The
novelty of our first result with respect to \cite{BenGG89} rests in
the uniformity of the constants with respect to small changes in the
frequencies $\omega$. Our scheme is reminiscent of that of Lochak in
his proof of Nekhoroshev's Theorem \cite{Loc92}: we use Dirichlet
approximation theorem in order to approximate the frequencies by
completely resonant ones and thus to reduce to a perturbation of a
completely resonant system. The error in the approximation by
Dirichlet theorem is controlled by a large parameter $\tQ$. Then we
apply the main theorem of \cite{\catene} which allows us to put the
system in resonant normal form up to a remainder which is
exponentially small with an effective small parameter.  Some work is
required in order to fit into the scheme of \cite{\catene}. Then one
gets a result in which there is an effective small parameter which
depends both on $\tQ$ and on $\epsilon$. So we choose $\tQ$ as a function
of $\epsilon$ in order to minimize the remainder, concluding the
proof.

The motivation for the second part of the paper rests in the remark
that if one is interested just in power law times, then the result can
be obtained by a normal form construction which is purely algebraic
(following the original ideas of Birkhoff). The only variant needed
with respect to the standard schemes is the one introduced by
\cite{\tedeschi}, namely to fix some threshold value $\alpha$ for the
small denominators and to consider as resonant all the monomials giving
rise to small denominators smaller than $\alpha$.  Then one can put
the system in resonant normal form (in the above sense) up to a
remainder of order $\epsilon^{-N}$. Finally, one has to prove that the
normal form admits an approximate integral of motion. We prove this last
fact using again Dirichlet theorem. We remark that the second result
holds also for Hamiltonians which are not analytic but only infinitely
differentiable.

\noindent
{\it Acknowledgments}. We thank Christian Lubich for pointing out a
mistake in the first version of the paper and for some comments that
led to considerable improvements of the paper.  This research was
founded by the Prin project 2010-2011 ``Teorie geometriche e
analitiche dei sistemi Hamiltoniani in dimensioni finite e infinite''.

%%%%%%%%%%%%%%%%%%%%%%%%%%%%%%%%%%%%%%%%%%%%%%%%%%%%%%%%%%%%%%%%%%%%%%%%
%
%               EXPONENTIALLY LONG TIMES
%
%%%%%%%%%%%%%%%%%%%%%%%%%%%%%%%%%%%%%%%%%%%%%%%%%%%%%%%%%%%%%%%%%%%%%%%%

\section{Exponentially long times}
\label{1}

In the phase space $\R^{2n}\oplus\R^{2d}\ni((p,q),(P,Q))$, endowed
with the usual euclidean norm, we consider a Hamiltonian system
of the form \eqref{sys.1}
%\begin{equation}
%\label{sys}
%H(p,q,P,Q)=h_\omega(p,q)+f(q,P,Q)\ ,
%\end{equation}
%where 
%\begin{equation}
%\label{homega}
%h_{\omega}(p,q):=\sum_{j=1}^{n}\frac{p_j^2+\omega_j^2q_j^2}{2}\ ,
%\end{equation}
where $H_0(P,Q,q)$ is \emph{analytic} in an open domain of $\R^{2d+n}$.

We first state the smoothness properties of $H_0$ in a precise form. For given
$E_0$ define the sublevel
\begin{equation}
\label{sub H}
\Sc_{E_0}:=\left\{ (P,Q)\in \R^{2d}\ :\ H_0(P,Q,0)\leq
E_0\right\}\ ,
\end{equation}
and the ball
\begin{equation}
\label{bcal}
\B_\rho:=
\left\{(p,q)\ :\ \norma{(p,q)}^2:=\sum_j\frac{p^2_j+q_j^2}{2}\leq
\rho^2\right\}\ .
\end{equation}
Remark that $\Sc_{E_0}$ needs not to be compact. Consider the
complexification of the phase space and denote by $B(\zeta,R)\subset
\C^{2n+2d}$ the closed ball of radius $R$ and center
$\zeta\equiv(p,q,P,Q)$.

We assume that {\it there exist positive $E_0^*,E^*, R^* $ such
  that, by defining
\begin{equation}
\label{dom}
\G:=\B_{3\sqrt{E^*}}\times \Sc_{3E_0^*}\ ,\quad
\G^*_{R^*}:=\bigcup_{\zeta\in\G}B(\zeta,R^*)\ ,
\end{equation}
 the function $H_0$ extends to a bounded analytic function on
 $\G^*_{R^*}$, namely to a function fulfilling
\begin{equation}
\label{e.sup_f}
\sup_{\G^*_{R^*}}\left|H_0(P,Q,q)\right| \leq C_{H_0}\ .
\end{equation}
} 

\begin{theorem}
\label{expo}
Under the above assumptions, there exist positive constants
$\epsilon_*, C_1,C_2$ such that, if $0<\epsilon<\epsilon_*$ and the initial
datum $(p^{0},q^{0},P^{0},Q^{0})$ fulfills
\begin{equation}
\label{dati}
H_0(P^0,Q^0,0)\leq E_0^{*}\ ,\quad h_\omega(p^0,q^0)\leq E^*\ ,
\end{equation}
then along the corresponding solution one has 
\begin{equation}
\label{estiexpo}
\left|h_\omega(t)-h_\omega(0)\right|<C_1\epsilon^{1/n}\ ,\quad
\text{for}\quad \left|t\right|\leq
C_2\exp\left(\frac{\epsilon_*}{\epsilon}\right)^{1/n} \ .
\end{equation}
The constants $\epsilon_*,$ $C_1$, $C_2$ depend only 
on $C_{H_0}$ and on $n$.
\end{theorem}

\begin{remark}
\label{unif}
The main point is that \emph{the constants do not depend on the
  frequencies and are thus uniform for all frequencies fulfilling
  \eqref{eps.1} with $\epsilon<\epsilon^*$.}
\end{remark}
%\begin{remark}
%\label{3r}
%The factor $3$ in the definition of the domains was inserted for
%simplicity. Everything could be reformulated eliminating such a factor
%and exploiting in a suitable way the fact that the domain $\G^*_{R^*}$
%striclty contains the domain $\G$. 
%\end{remark}

\proof First we remark that by Cauchy inequality for analytic
functions one has that the quantities
\begin{equation}
\label{grad}
\left|\frac{\partial H_0}{\partial
  q_j}\right|\ ,\quad\left|\frac{\partial H_0}{\partial
  Q_l}\right|\ ,\quad\left|\frac{\partial H_0}{\partial P_l}\right|\
\end{equation}
are bounded on any domain contained in $\G^*_{R^*} $, hence the same holds true for the Hamiltonian vector field
$X_{H_0}$.

To be definite we assume
\begin{equation}
\label{e.omega1}
\min\{\omega_j\}=\omega_1=\frac1\epsilon\ .
\end{equation}
According to Dirichlet theorem, for any $\tQ>1$ there exist
integers $\tq\leq \tQ$ and $\{\tp_j\}_{j=2}^n$ s.t.
\begin{equation}
\label{diri}
\left|\frac{\omega_j}{\omega_1}-\frac{\tp_j}{\tq}\right|\leq
\frac{1}{\tq\tQ^{1/(n-1)}} \ ,\quad j=2,...,n\ .
\end{equation}
The value of $\tQ$ will be fixed later on as a function of $\epsilon$.

Define a new vector of resonant frequencies $\tilde\omega$
\begin{equation}
\label{e.om.tilde}
\tilde \omega_1:=\omega_1\ ,\qquad \tilde\omega_j:=\omega_1
\frac{\tp_j}{\tq}=\frac{\tp_j}{\epsilon \tq}\ ,
\end{equation}
and  
\begin{equation}
\label{acche}
h_{\tilde \omega}:=\sum_{j=1}^{n}\frac{p_j^2+\tilde
  \omega_j^2q_j^2}{2}\ ,\quad  h_1(q):=\frac{1}{2}
\sum_{j=1}^{n}(\omega_j^2-\tilde\omega_j^2)q_j^2\ ,\quad f:=H_0+h_1
\end{equation}
so that the Hamiltonian takes the form 
\begin{equation}
\label{tilde}
H=h_{\tilde\omega}+f\ ,
\end{equation}
as required in \cite{BamG93}. Then \eqref{diri} becomes
\begin{equation}
\label{e.diri.tilde}
\left|\frac{\omega_j - \tilde\omega_j}{\omega_1}\right|\leq
\frac{1}{\tq\tQ^{1/(n-1)}} \ ,\quad j=2,...,n\ .
\end{equation}
Furthermore the flow generated by $h_{\tilde \omega} $ is periodic
with frequency $\omega:=1/\epsilon \tq$.

We redefine the norms and the domains in order to fit the scheme by
\cite{BamG93}. So we put
\begin{align}
\label{nuovenorme}
\normati{\zeta}^2 &\equiv\normati{(p,q,P,Q)}^2 =
\sum_{j=1}^{n}\frac{|p_j|^2+\tilde\omega_j^2|q_j|^2}{2}+
\sum_{l=1}^{d}
\frac{|P_l|^2+|Q_l|^2}{2}\equiv\\ &\equiv\normati{(p,q)}^2+\norma{(P,Q)}^2
\nonumber\\
\label{nuG}
\tilde \G &:=  \tilde\B_{3\sqrt{E^*}}  \times \Sc_{3E^*_0}\ ,\quad
\tilde\G_{\tilde R}:=\bigcup_{\zeta\in\tilde \G}\tilde B(\zeta,\tilde
R)\ ,
\end{align}
where $\tilde \B$ and $\tilde B$ are the closed ball in
the norm \eqref{nuovenorme}.

The relation with the old norms and domains is easily obtained: 
the new norm \eqref{nuovenorme} is stronger than the euclidean one:
$$
\norma{q}\leq \epsilon \normati{q}\ ,
$$ so, provided $\epsilon$ is small enough and $\tQ$ large enough,
choosing  $\tilde
R:=R^*/2$, one
has $\tilde \G_{\tilde R}\subset \G^{*}_{R^*}$ (strictly and with some finite distance between the
boundaries).

We have now to compute the constants involved in the statement of
Theorem 4.1 of \cite{BamG93}, namely
\begin{align*}
\omega_f:=\frac{1}{\tilde R}\sup_{\zeta\in\tilde \G_{\tilde
    R}}\normati{X_f(\zeta)}\leq \frac{1}{\tilde R}\left[\sup_{\zeta\in\tilde \G_{\tilde
    R}}\normati{X_{H_0}(\zeta)}+\sup_{\zeta\in\tilde \G_{\tilde
    R}}\normati{X_{h_1}(\zeta)}\right]\ .
\end{align*}
Using \eqref{grad} and \eqref{e.sup_f} one immediately sees that the
supremum of $X_{H_0} $ is independent of $\epsilon$ and of ${\tQ}$. In
order to compute the supremum of $X_{h_1}$ recall
\eqref{e.diri.tilde}, and remark that
\begin{align*}
\frac{\omega_j}{\tilde
  \omega_j}=1+\frac{\omega_j-\tilde\omega_j}{\tilde \omega_j}\ ,\quad
\left| \frac{\tilde\omega_j-\omega_j}{ \tilde \omega_j}\right|
=\frac{\left|\frac{\tilde\omega_j-\omega_j}{\tilde \omega_1}\right|}
{\big|\frac{\omega_j}{\omega_1}-\frac{\omega_j-\tilde\omega_j}{\omega_1}\big|}
\leq\frac{1/\tq\tQ^{1/(n-1)}}{1-\frac{1 }{\tq\tQ^{1/(n-1)}}} \leq
\frac{2}{\tq\tQ^{1/(n-1)}}
\end{align*}
provided $\tq\tQ^{1/(n-1)}>2$, from which 
\begin{equation}
\label{e.est.1}
\left|\frac{\tilde\omega_j+\omega_j}{\tilde\omega_j}\right|\leq
3\ ,\qquad\qquad
\left|\frac{\tilde\omega_j^2-\omega_j^2}{\tilde\omega_j\omega_1}\right|
\leq \frac{6}{\tq\tQ^{1/(n-1)}}\ .
\end{equation}
Thus using \eqref{nuovenorme} and \eqref{e.est.1} the field $X_{
  h_1}$ admits the upper bound
\begin{align*}
\normati{X_{h_1}(\zeta)}^2 &=
\sum_{j=1}^{n}(\omega_j^2-\tilde\omega_j^2)^2q_j^2= \sum_{j=1}^{n}
\frac{(\omega_j^2-\tilde\omega_j^2)^2}{\tilde\omega_j^2\omega_1^2}
\omega_1^2\tilde\omega_j^2 q_j^2\leq \\ &\leq
\frac{1}{\epsilon^2}\quadr{\sup_{j=1,...,n}
  \Big|\frac{\omega_j^2-\tilde\omega_j^2}{\tilde\omega_j\omega_1}\Big|^2}
\sum_{j=1}^{n}\tilde\omega_j^2 q_j^2 \leq
\left(\frac{6}{\tq\tQ^{1/(n-1)}}\right)^2
\frac{2}{\epsilon^2}\normati{(p,q)}^2\ ,
 \end{align*}
which gives
\begin{displaymath}
\sup_{\zeta\in\tilde \G_{\tilde R}}\normati{X_{h_1}(\zeta)} \leq
\left(\frac{6\sqrt2}{\tq\tQ^{1/(n-1)}}\right)
\frac{\sqrt{9E^*+{\tilde R}^2}}{\epsilon}\ .
\end{displaymath}
So one can put
\begin{equation}
\label{e.omega.0}
\omega_f \leq C\left[1 + \frac{1}{\epsilon \tq\tQ^{1/(n-1)}}\right]\ ,
\end{equation}
and the small parameter $\mu$ of Theorem 4.1 turns out to
be\footnote{Recall that in Dirichlet Theorem $\tq\leq \tQ$.}
\begin{equation}
\label{mu}
\mu:= C{\frac{\omega_f}{\omega}} \leq C\epsilon
\tq\left(1+ \frac{1}{\epsilon \tq\tQ^{1/(n-1)}} \right)\leq
C_1\left(\epsilon \tQ+ \frac{1}{\tQ^{1/(n-1)}}\right)\ .
\end{equation}
Following \cite{BamG93} p.~604, we choose
$\tQ^{1/(n-1)}=\epsilon^{-1/n}$, so that $\epsilon \tQ = \epsilon^{1/n}$
and we can choose $\mu= C_2\epsilon^{1/n}$. Defining
$\epsilon_*:=C_2^{-n}$ and computing the other constants in Theorem
4.1 and its corollaries one gets the thesis. \qed

%%%%%%%%%%%%%%%%%%%%%%%%%%%%%%%%%%%%%%%%%%%%%%%%%%%%%%%%%%%%%%%%%%%%%%%%
%
%               POWER LAW
%
%%%%%%%%%%%%%%%%%%%%%%%%%%%%%%%%%%%%%%%%%%%%%%%%%%%%%%%%%%%%%%%%%%%%%%%%

\section{Power law times}
\label{2} 

\subsection{Statement}

The aim of this section is to give an easy proof of a simplified
result, in which the control of the energy of high frequency
oscillators is obtained only for time scales of order $\epsilon^{-N}$
with an arbitrary $N$. We remark that for the present result
$C^\infty$ smoothness is enough. Precisely Theorem \ref{main} below is
true under the assumption that there exists an interval of values of
$E_0$, $\rho$ s.t., for any $k$ the $C^k$ norm of $H_0$ is bounded in
$\B_{3\rho}\times \Sc_{3E_0}$. Of course, if one fixes a value of $N$ then
finite smoothness is also enough.

\begin{theorem}
\label{main}
Fix a positive (small) $b$, then, for any positive (large) $N$,
there exists a positive constant ${\epsilon}_*(N,b)$, such that, if
${\epsilon}<{\epsilon}_*$, and the initial datum fulfills
\begin{equation}
\label{ini.dat}
E:=h_{\omega}(p,q)<E^*\ ,\qquad H_0(P,Q,0)<E_0^*\ ,
\end{equation}  
then one has
\begin{equation}
\label{esti}
| h_{\omega}(t)-h_\omega(0)| \leq E b\ \qquad\text{for}\ |t|\leq
{\epsilon}^{-N}\ .
\end{equation}
\end{theorem} 

\begin{remark}
\label{finitedim}
The constant ${\epsilon}_*$ strongly depends on the dimension $n$ of
the fast system, so the result does not extend to infinite dimensional
systems.
\end{remark}
\begin{remark}
\label{sub}
In the present statement the change of the energy of the high
frequency system is controlled by the parameter $b$, which is arbitrarily
small, but independent of $\epsilon$. On the contrary, in the paper
\cite{\tedeschi} one has $b\sim\epsilon^{3/4}$.
\end{remark}

%%%%%%%%%%%%%%%%%%%%%%%%%%%%%%%%%%%%%%%%%%%%%%%%%%%%%%%%%%%%%%%%%%%%%%%%
%
%               PROOF
%
%%%%%%%%%%%%%%%%%%%%%%%%%%%%%%%%%%%%%%%%%%%%%%%%%%%%%%%%%%%%%%%%%%%%%%%%

\subsection{Proof}
\label{proof}

We start by preparing the Hamiltonian, subsequently we introduce the
kind of expansion needed for the further developments. Then we
prove an approximation lemma for the frequencies and finally we
prove the normal form lemma that we will use to get Theorem \ref{main}.

First we scale the variables, the frequencies and the time in a
suitable way (see also Sections 2 and 4 of \cite{BenGG89}). Together,
we introduce the standard complex variables usually needed in order to
develop perturbation theory. As in sect. \ref{1} we assume
$\omega_1=\min{\omega_j}$. 

Thus define
\begin{equation}
\label{scaling}
\nu_j:={\epsilon}\omega_j\ ,\qquad
p_j=\sqrt\frac{\nu_j}{{2\epsilon}}\frac{\xi_j+\eta_j}{i} \ ,\qquad 
q_j=\sqrt{\frac{\epsilon}{2\nu_j}} (\xi_j-\eta_j) \ ,
\end{equation}
(in particular one has $\nu_1=1$) so that, by rescaling time to $t' :=
\epsilon t$, the Hamiltonian of the system (still denoted by $H$)
takes the form
\begin{equation}
\label{hres}
H=\sum_{j=1}^{n}\nu_j\xi_j\eta_j + {\epsilon}
H_0(P,Q,q(\xi,\eta))\ .
\end{equation}
For the new fast variables $(\eta,\xi)$ we will use the norm
\begin{equation}
\label{norm.ph}
\norma{(\xi,\eta)}^2 =
\sum_{j=1}^{n}{\nu_j}\tond{|\xi_j|^2+|\eta_j|^2} 
\ .
\end{equation}
which corresponds to the rescaled energy norm in the original $(p,q)$
variables.

If we define
\begin{equation}
\label{def.rhoE}
\rho_E:=\sqrt {E}\ ,
\end{equation}
then $\sum_j(p_j^2+\omega_j^2q_j^2)\leq E$ implies $h_\nu\leq
\epsilon \rho_E^2$, which means $(\xi,\eta)\in
\B_{\rho_E\sqrt{\epsilon}}$. Hence, since the variables $\xi,\eta$
have size of order $\sqrt{\epsilon}$, we have to consider an expansion
of the nonlinear terms in both $\sqrt{\epsilon}$ and in $\xi,\eta$. In
other words, the scaling \eqref{scaling} introduces two different
dependencies on $\sqrt{\epsilon}$ in the Hamiltonian: an implicit one,
of size $\mathcal{O}(\sqrt{\epsilon })$ in the scaled variables, and
an explicit one in the coefficient in front of any monomial depending
on $(\xi,\eta)$, due to dependence on $q$ only.

As anticipated above the main step of the proof consists in putting
the system in normal form. We now specify in a precise way what we
mean by normal form.

\begin{definition}
\label{normal}
Given $\alpha>0$, a monomial
$\xi^l\eta^m\equiv\xi_1^{l_1}....\xi_n^{l_n}\eta_1^{m_1}....\eta_n^{m_n}$
is said to be in $\alpha$-normal form if
\begin{equation}
\label{nor.def}
\left|\nu\cdot(l-m)\right|\leq \alpha\ .
\end{equation}
\end{definition}

We are now going to prove that, if $\alpha$ is small enough, then
there exists a {\it non vanishing} vector $\tilde \nu$ such that
\begin{equation}
\label{tilde.1}
h_{\tilde \nu}(\xi,\eta):=\sum_{j=1}^{n}\tilde \nu_j\xi_j\eta_j
\end{equation}
Poisson commutes with all the monomials in normal form.

%%%%%%%%%%%%%%%%%%%%%%%%%%%%%%%%%%%%%%%%%%%%%%%%%%%%%%%%%%%%%%%%
%
%          APPROSSIMAZIONE RAZIONALE DELLE FREQ
%
%%%%%%%%%%%%%%%%%%%%%%%%%%%%%%%%%%%%%%%%%%%%%%%%%%%%%%%%%%%%%%%%

\begin{lemma}
\label{cor.tilde}
Fix $N>0$, then there exists a non negative sequence
$\{\alpha_i\}_{i\geq 1} $, with $\lim_{i\to\infty}\alpha_i=0$, such
that, for every frequency vector $\nu$ there exists a new frequency
vector $\tilde \nu$, depending on $\alpha_i$, which fulfills
\begin{align}
\label{f.f}
\sup_{j=1,\ldots,n}|\tilde\nu_j-\nu_j|\leq \frac{\alpha_i}{N}\ , \\
\label{f.f.1}
\left\{h_{\tilde \nu};\xi^l\eta^m\right\}=0\ ,
\end{align}
for all monomials $\xi^l\eta^m$ in $\alpha_i$ normal form satisfying
$|l|+|m|\leq N$.
\end{lemma}

\proof We use again Dirichlet theorem. The form we choose is the one
according to which,
for any $\nu\in\R^{n-1}$ the inequalities
\begin{equation}
\label{diri.1}
\left|\nu_j-\frac{\tp_j}{\tq}\right|\leq\frac{1}{\tq^{1+1/(n-1)}} \ , \quad
j=2,...,n
\end{equation}
have infinitely many solutions $\tq\in\Na$, $\tp_j(\tq)\in\Z$. In particular
the $\tq$'s form a diverging sequence $\tq_i$. We identify the sequence
$\tq_i$ with the corresponding value of $\tq$ (instead of using
$i$). Define $\tilde\nu _j:=\tp_j/\tq$, $j=2,...,n-1$, $\tilde
\nu_1:=\nu_1=1$ and $\alpha_{\tq}(N,n):=N/\tq^{1+1/(n-1)}$. We are now going to
prove that $\tilde \nu\cdot k\not=0$ with $|k|\leq N$ implies
$|\nu\cdot k|>\alpha_{\tq}$. First remark that $\tilde\nu\cdot k\not=0$
implies $|\tilde\nu\cdot k|\geq 1/\tq $ (since $\tilde \nu_j$ are
rationals), so that one has
\begin{align*}
|\nu\cdot k|&\geq |\tilde\nu\cdot k|-\left|(\nu-\tilde\nu)\cdot
k\right|\geq \frac{1}{\tq}-|\tilde\nu-\nu||k| \geq\\
&\geq \frac{1}{\tq}-\frac{N}{\tq^{1+1/(n-1)}}=\frac{1}{\tq}-\alpha_{\tq}
=\alpha_{\tq}\left[\frac{1}{(N^{n-1}\alpha_{\tq})^{\frac{1}{n}}}-1\right]\ ,
\end{align*}
but, provided $\alpha_{\tq}$ is small enough with respect to $N^{n-1}$ the
square bracket is bigger than 1 and the thesis follows.  \qed

\bigskip

We fix now once for all $\alpha$ as 
\begin{equation}
\label{e.def.alpha}
\alpha \equiv\alpha_i\leq
\frac{bN}{21}\ .
\end{equation}

In the following we are going to construct a canonical transformation
which puts the Hamiltonian in normal form up to order
$(\sqrt{\epsilon})^{2N} = \epsilon^N$.  We first introduce the class
of polynomials that we will meet in the construction and the degree
that we will assign to each of them.

\begin{definition}
\label{ps}
Let $\U\subset\Sc_{3E^*_0}$ be an open domain. For $s\geq 0$, the space
$\Ph_s\equiv\Ph_s(\U)$ is the space of the linear combinations, with
coefficients in $C^\infty(\U)$, of the
monomials of the form
\begin{equation}
\label{mono}
(\sqrt{\epsilon})^{a+2}\xi^l\eta^m\ ,
\end{equation}
where the indexes fulfill the conditions
\begin{equation}
\label{indi}
a+|l|+|m|=s\ ,\quad a\geq |l|+|m|\ .
\end{equation}
If $g\in\Ph_s$, then the index $s$ will be called the {\sl order} of
the polynomial $g$.
\end{definition}

In the following, when not needed, we will not specify the domain $\U$.
It is immediate to verify the following Lemma
\begin{lemma}
\label{l.poisson}
Let $g_1\in\Ph_{s_1}$ and $g_2\in\Ph_{s_2}$, then
\begin{equation}
\label{poisson}
\left\{g_1;g_2\right\} \in \Ph_{s_1+s_2}\oplus \Ph_{s_1+s_2+2}\ .
\end{equation}
\end{lemma}

\proof Indeed
\begin{displaymath}
\left\{g_1;g_2\right\} = \left\{g_1;g_2\right\}_{P,Q} +
\left\{g_1;g_2\right\}_{\xi,\eta}\ ;
\end{displaymath}
the first term at r.h.s. belongs to $\Ph_{s_1+s_2+2}$ and the second
one belongs to $\Ph_{s_1+s_2}$. \qed

\begin{remark}
\label{rm.poi} 
Consider
$\left\{h_\nu;g\right\}$, with $g\in\Ph_{s}$ and $s\geq 1$. In this
case it, due to the lack of a prefactor $\epsilon$ in front of
$h_\nu$, is immediate to verify that $\left\{h_\nu;g\right\}\in
\Ph_s$.
\end{remark}

\begin{remark}
\label{2r}
Moreover, it is useful to stress that both in case of Lemma
\ref{l.poisson} and in the case of $\left\{h_\nu;g\right\}$, the
parity of the space $\Ph_s$ is preserved by the Poisson brackets. Due
to the structure of the perturbation $\epsilon H_0$, we will deal only
with even parity spaces $\Ph_{2s}$.
\end{remark}

\bigskip
It is useful to extend the definition to functions of
$\xi,\eta,\sqrt\epsilon$ of class $C^\infty$
and to introduce the space of the functions that will play the role of
remainders.

\begin{definition}
\label{fu}
Let $F((P,Q),(\xi,\eta),\sqrt{\epsilon})$, $F\in C^{\infty}(\U\times
\B_{\rho}\times \B_{\sqrt{\epsilon^\sharp}})$ for some positive
$\rho,\sqrt{\epsilon^\sharp}$. We say that $F\in\bar \Ph(\U)$ if each
of its Taylor polynomials in $\xi,\eta,\sqrt{\epsilon}$ belongs to
some of the spaces $\Ph_s(\U)$.
\end{definition}

Given a function $F\in\bar \Ph$ we can define the projector
$\Pi_s$ which extracts from $F$ its component in $\Ph_s$.

\begin{definition}
\label{res}
A function $F\in \bar\Ph(\U)$ will be said to belong to $\resto_r(\U)$
if one has $\Pi_sF=0$, $\forall s\leq r$.
\end{definition}

\begin{remark}
For any $N$ one can expand $H_0$ in Taylor series in the variables
$\xi,\eta$ at order $N$, getting
$$
\epsilon H_0=\sum_{s=0}^{N}f_s+R^{(N)}\ ,\quad f_s(P,Q,\xi,\eta)= 
{\epsilon}\sum_{|l|+|m|=s}a_{lm}(P,Q)\xi^l\eta^m\epsilon^{s/2} 
$$
and $R^{(N)}$ having a zero of order $N+1$ in the variables
$\xi,\eta$. Thus one has $f_s\in\Ph_{2s}$ and
$R^{(N)}\in\resto_{2N+1}$ (see Remark \ref{2r}).
\end{remark}

\bigskip
\begin{remark}
\label{resto}
Let $F\in \resto_{2s+1}(\U)$, with $\U\subset \Sc_{3E^*_0}$, then one
has $\sup_{\U\times B_{\rho\sqrt{\epsilon}}}|F|\leq
C\sqrt{\epsilon}^{2s+3}$. The constant depends in particular on $\U$
and on $\rho$. Similar inequalities hold for the derivatives of $F$.
\end{remark}

The normalizing transformation will be constructed using the Lie
transform $\phi_{\chi}$, namely the time one flow of an auxiliary
Hamiltonian $\chi\in\Ph_{2r}$ with $r\geq 1$. The main properties of the
Lie transform are summarized in the next lemma.
\begin{lemma}
\label{lie}
Let $\Sc_{3E^*_0}\supset\U_1\supset\U_2\supset\U_3\supset\Sc_{2E^*_0}$
be open sets (the inclusion must be strict) and let
$3\rho_E>\rho_1>\rho_2>\rho_3>2\rho_E$ be positive parameters. Let
$\chi\in\Ph_s(\U_1)$ with $s\geq 1$. Then there exists
${\epsilon}_{\sharp}$, such that, if ${\epsilon}<{\epsilon}_{\sharp}$,
then one has
\begin{equation}
\label{defo}
\U_1\times B_{\rho_1\sqrt{\epsilon}}\supset \phi_{\chi}(\U_2\times
B_{\rho_2\sqrt{\epsilon}}) \supset \U_3\times B_{\rho_3\sqrt{\epsilon}}\ .
\end{equation} 
The constant ${\epsilon}_\sharp$ depends only on the above sets
$\U_i$, on $\rho_i$ and on the norm ${\cal C}^1(\U_1)$ of the coefficients of
the development of $\chi$ in $\xi,\eta,\sqrt\epsilon$.

Let $F\in \Ph_r$, then one has
\begin{equation}
\label{lie.tr}
\left[F\circ\phi_{\chi}-F\right]\in \resto_{s+r}
\end{equation}
and 
\begin{equation}
\label{lie.tr1}
\left[h_\nu\circ\phi_{\chi}-\left(h_\nu+\left\{\chi;h_\nu\right\}
  \right)\right]\in \resto_{2s-1}\ .
\end{equation}
\end{lemma}
\proof The statement on the existence of the flow and the way it
transforms open domains immediately follows from the standard theory
of existence and uniqueness of ODEs.

To get \eqref{lie.tr} and \eqref{lie.tr1} one uses   
$$
\frac{\di}{\di t}F\circ \phi^t_{\chi}=\poisson \chi F\circ
\phi^t_\chi\ ,
$$
from which, 
\begin{equation}
\label{restotaylor}
F\circ \phi_\chi=F+\left\{
\chi;F\right\}+\int_0^1(1-s)\left\{\chi;\left\{\chi;F\right\}\right\}
\circ\phi_{\chi^s_{\chi}}\, \di s\ ,
\end{equation}
which holds both for the case of the function $F$ of the statement and
for the function $h_\nu$. Then using Lemma \ref{l.poisson}, the fact
that $\left\{\chi;h_\nu\right\}\in\Ph_s(\U_1)$ (see Remark
\ref{rm.poi}) and standard estimates the thesis follows.  \qed

We are now ready to state and prove the iterative lemma which yields
the existence of the normal form.

\begin{lemma}
\label{norm.form.lemma}
There exists a sequence of domains
$\Sc_{3E^*_0}\supset\U_0\supset\U_1\supset...\supset
\U_{N+1}\supset\Sc_{2E^*_0}$ and a sequence of positive parameters
$3\rho_E>\rho_0>\rho_1>...>\rho_{N+1}>2\rho_E$ with the following
property: for any $1 \leq r\leq N$ there exists a positive
${\epsilon}_r$, such that, if ${\epsilon}<{\epsilon}_r$ then there
exists a canonical transformation $T^{(r)}:\U_r\times
B_{\rho_r\sqrt{\epsilon}}\to \U_0\times B_{\rho_0\sqrt{\epsilon}}$,
$T^{(r)}( \U_r\times B_{\rho_r\sqrt{\epsilon}} )\supset \U_{r+1}\times
B_{\rho_{r+1}\sqrt{\epsilon}}$ such that $H\circ T^{(r)}$ is in normal
form at order $2r$, namely $\forall l\leq 2r$ the polynomial
$\Pi_{l}\left[H\circ T^{(r)}\right]$ is in normal form. One also has
\begin{equation}
\label{defo1}
\left[h_{\tilde \nu}\circ T^{(r)}-h_{\tilde
    \nu}\right]\in\resto_1\ ,\quad \left[\epsilon H_0\circ
  T^{(r)}-\epsilon H_0\right]\in\resto_3\ .
\end{equation}
The sets $\U_j$, as well as the parameters $\rho_j$, the $C^k$ norm of
$T^{(r)}$, and the quantity ${\epsilon}_r$, depend on the frequency
only through the parameter $\alpha$. Finally the transformed
Hamiltonian contains only terms of even order (in the sense of
definition \ref{ps}).
\end{lemma}

\proof The proof follows the standard proof of Birkhoff normal form
theorem. The theorem is true for $r={0}$. We assume it for $r$ and
prove it for $r+1$. We construct the transformation increasing by one
the order of the non normalized part of the Hamiltonian as the Lie
transform generated by a function $\chi_{r+1}\in \Ph_{2(r+1)}$. First
remark that, according to \eqref{lie.tr}, \eqref{lie.tr1} the
transformed Hamiltonian is automatically in normal form at order
$2r$. We are going to choose $\chi_{r+1}$ in such a way that
$$
\Pi_{2(r+1)}\left[H\circ T^{(r)}\circ \phi_{\chi_{r+1}} \right] \equiv
\left\{\chi_{r+1};h_\nu \right\}+ \Pi_{2(r+1)}\left[H\circ T^{(r)} \right]
$$
is in normal form too. To this end write
$$
\Pi_{2(r+1)}\left[H\circ T^{(r)}
  \right]=\sum_{a,l,m}P_{a,l,m}(P,Q)\sqrt{\epsilon}^{a+2} \xi^{l}\eta^m\ ,
$$ where the indexes fulfill the limitations $a+|l|+|m|=2(r+1)$ and $a
\geq |l|+|m|$, so that, in particular $|l|+|m|\leq r+1\leq N$. Define
now
\begin{equation}
\label{chi}
\chi_{r+1}:=\sum_{(l,m)\in NR,\ a}\frac{P_{a,l,m}(P,Q)}{\im\nu\cdot
  (l-m)}\sqrt{\epsilon}^{a+2} \xi^{l}\eta^m
\end{equation}
where the nonresonant set $NR$ is defined by
\begin{equation}
\label{nres}
NR:=\left\{ (l,m)\ :\ \left|\nu\cdot (l-m)\right|>\alpha\right\}\ .
\end{equation}
Then, the $C^k$ norm of $\chi_{r+1}$ is controlled by the $C^k$ norm
of $\Pi_{2(r+1)}\left[H\circ T^{(r)} \right]$ divided by $\alpha$, and
therefore the statement on the $C^k$ norm of the transformation
holds. The statement on the domain of definition of the transformation
follows from $\U_{r+2}\subset
\phi_{\chi_{r+1}}(\U_{r+1})\subset\U_{r}$ which is consequence of
Lemma \ref{lie}. The same is true for \eqref{defo1} (which at leading
order follows directly from Lemma \ref{l.poisson}) and the statement
on the dependence of the parameters on the frequency.  \qed

\bigskip

%%%%%%%%%%%%%%%%%%%%%%%%%%%%%%%%%%%%%%%%%%%%%%%%%%%%%%%%%%%%%%%%
%
%          FINE DELLA DIMOSTRAZIONE
%
%%%%%%%%%%%%%%%%%%%%%%%%%%%%%%%%%%%%%%%%%%%%%%%%%%%%%%%%%%%%%%%%

\noindent{\bf End of the proof of Theorem \ref{main}.} Consider
$T^{(N)}$ and denote the new variables by $(P',Q',\xi',\eta')$, namely
$(P,Q,\xi,\eta)=T^{(N)}(P',Q',\xi',\eta')$ and by $h'_{\tilde \nu}:=
\sum \tilde\nu_j\xi'_j\eta'_j$. Compute
\begin{align}
\label{muov}
\left| h_\nu(t)-h_\nu(0)\right|&\leq \left| h_{\nu}(t)-h_{\tilde
  \nu}(t)\right|+\left| h_{\tilde \nu}(t)-h'_{\tilde \nu}(t)\right| 
\\
\label{muov.1}
&+ \left| h'_{\tilde \nu}(t)-h'_{\tilde \nu}(0)\right|+\left| h_{\tilde
  \nu}'(0)-h_{\tilde \nu}(0)\right| +\left| h_{\tilde
  \nu}(0)-h_{\nu}(0)\right| \ .   
\end{align} 
Assume for a while that $\norma{(\xi'(t),\eta'(t))}\leq
2\rho_E\sqrt{\epsilon}$ for $|t|\leq \epsilon^{-N}$, then from Lemma
\ref{norm.form.lemma} $\norma{(\xi(t),\eta(t))}\leq
3\rho_E\sqrt{\epsilon}$ and one can use \eqref{f.f} and \eqref{defo1}
to estimate the different terms of \eqref{muov} and \eqref{muov.1} by
\begin{equation}
\label{muov.2}
4\rho_E^2{\epsilon}\frac{\alpha }{N} +
C{\epsilon}^2+\left|\left\{h_{\tilde\nu}, H\circ
T^{(N)}\right\}\right| |t|+C{\epsilon}^2 +
9\rho_E^2{\epsilon}\frac{\alpha }{N} \ .
\end{equation}
Indeed from $h_{\tilde\nu} - h'_{\tilde \nu}\in\resto_1$ it follows
immediately $| h_{\tilde \nu}(t)-h'_{\tilde \nu}(t)| <
C\epsilon^2$. On the other hand, one has to recall that $h_\nu$ is the
norm (see \eqref{norm.ph}) and that $h_{\tilde\nu}$ is close to
$h_\nu$ because of \eqref{f.f}
\begin{displaymath}
\left|h_\nu-h_{\tilde\nu}\right| \leq
\tond{\sup_{j=1,\ldots,n}|\nu_j-\tilde\nu_j|}\sum_j\tond{|\xi_j|^2+|\eta_j|^2}
\leq \frac{\alpha }{N} \|(\xi,\eta)\|^2\ .
\end{displaymath}

Now, since $H\circ T^{(N)}$ is in normal form, one has
$\left\{h_{\tilde\nu}, H\circ T^{(N)}\right\}\in\resto_{2N+1}$ which in
turn implies $\left|\left\{h_{\tilde\nu}, H\circ
T^{(N)}\right\}\right|\leq C{\epsilon}^{N+2}$ and therefore, for the
considered times the third term is smaller than $C {\epsilon}^2$.

Take now ${\epsilon}$ so small that the sum of the second, the
third and the forth term of \eqref{muov},\eqref{muov.1} does not
exceed $\alpha \rho_E^2 {\epsilon}/N$, then going back to the original
variables and recalling that, from \eqref{e.def.alpha},
$\alpha\leq Nb/21 $ the estimate \eqref{esti} follows.

We still have to prove that for $|t|\leq \epsilon^{-N}$ all the variables
are in the domain of validity of the normal form. Concerning the fast
variables this is a consequence of an argument similar to that of
Lyapunov's theorem which gives
$$ h'_{\tilde \nu}(t)\leq h'_{\tilde \nu}(0)+\left| h'_{\tilde
  \nu}(t)-h'_{\tilde \nu}(0)\right| \leq \rho_E^2
{\epsilon}(1+\alpha)+C{\epsilon}^2 \leq
2\rho_E^2 {\epsilon}\ .
$$ Concerning the variables $(P',Q')$ we exploit the conservation of
the Hamiltonian. To this end denote $\hat h(P,Q):=H_0(P,Q,0)$ and
$H_P:=H_0-\hat h$, and remark that $|H_P|<C\epsilon$, so that one has
(in the $(P,Q)$ variables)
$$
\hat h(t)=\hat h(0)+h_\omega(0)-h_{\omega}(t)+H_P(0)-H_P(t)
$$
so that, recalling the second of \eqref{e.def.alpha}, one has 
$$ \hat h(t)\leq E^*_0+Eb+C {\epsilon}<\frac{3}{2}E^*_0\ ,
$$ provide $b$ and $\epsilon$ are small enough.  It follows that $\hat
h'(t)\leq 2E^*_0$ on the considered time scale. The result then holds
in the rescaled time. To get the result in the physical time, just
repeat the whole argument with $N+1$ in place of $N$.\qed

%\bibliographystyle{amsalpha}
%\bibliographystyle{plain}
%\bibliography{database}

\providecommand{\bysame}{\leavevmode\hbox to3em{\hrulefill}\thinspace}
\providecommand{\MR}{\relax\ifhmode\unskip\space\fi MR }
% \MRhref is called by the amsart/book/proc definition of \MR.
\providecommand{\MRhref}[2]{%
  \href{http://www.ams.org/mathscinet-getitem?mr=#1}{#2}
}
\providecommand{\href}[2]{#2}

\end{document}